\documentclass{article}

\usepackage{amsmath}
\usepackage[francais]{babel}
\usepackage{amssymb}
\usepackage[applemac]{inputenc}

\def\dim{\mathop{\rm dim}}
\def\car{\mathop{\rm char}}
\def\s{\mathop{\rm ss}}
\def\GL{\mathop{\rm GL}}
\def\SL{\mathop{\rm SL}}

\def\End{\mathop{\rm End}}
\def\Gal{\mathop{\rm Gal}}
\def\Hom{\mathop{\rm Hom}}
\def\Ker{\mathop{\rm Ker}}
\def\I{\hbox{\rm Im}}
\def\Z{\bf Z}
\def\Q{\bf Q}

\def\e{\epsilon}
\def\+{\big{\mathop{\boldsymbol{+}}}}

\def\n{\noindent}
\def\m{\frak{m}}

\begin{document}

\begin{center}{\large On the mod $p$ reduction of orthogonal representations} 

\end{center}

\medskip

\begin{center}  Jean-Pierre Serre

\end{center}  
\medskip

\begin{center}  {\it to the memory of Bertram Kostant}

\end{center}  

\bigskip  
\
{\bf Introduction}

\medskip
\n   The following fact is probably known :  
  
  \smallskip
  
 (F) \ \     {\it Let  $G$  be a finite group and let $\rho$ be a linear representation of  $G$  in characteristic $0$ which is orthogonal} (resp. symplectic). {\it Let  $p$  be a prime number $\neq 2$.  Then the reduction mod $p$ of $\rho$ is orthogonal} (resp. symplectic).  
     
     \smallskip
     
 \n {\small [Recall that a linear representation is said to be orthogonal (resp. symplectic) if it fixes a nondegenerate symmetric (resp. alternating) bilinear form. As for the precise meaning of ``reduction mod $p$'', see below.]}
 
 \smallskip 
     In what follows, we give a detailed proof of (F), and we extend it to linear representations of algebras with involution, including the group algebra of an infinite group. 
     
     \medskip
     
     Let us state the theorem explicitly, in the group algebra case:

     \medskip
     Let  $R$ be a discrete valuation ring, with ring of fractions  $K$  and residue field $k=R/\pi R$,
     where  $\pi$  is a uniformizer. Let $V$ be a finite dimensional $K$-vector space. Let $G$ be a group and let $\rho : G \to \GL(V)$ be a homomorphism. Assume that there exists a lattice
      $L$ of  $V$  which is $G$-stable  (such a lattice always exists when $G$ is finite).
   
   \smallskip
    {
\n   \small[Recall that a {\it lattice} is a free $R$-submodule $L$ of $V$ such that $K.L = V$, cf. §1.1.]}

\smallskip
  Let  $L$ be a $G$-stable lattice of $V$. The $k$-vector space $V_L = L/\pi L$  is a  $k[G]$-module. The structure of this module may depend on the choice of $L$.  However, by a theorem of Brauer-Nesbitt 
(see §2), its {\it
semisimplification}  $V_L^{\s} $ is well-defined, up to isomorphism (cf. §2.2). It is $V_L^{\s} $ that we call ``the reduction mod $\pi $''  of $V$, and we denote it by $V_k$. The precise form of (F) is :

\smallskip

{\bf Theorem A}. {\it Assume that there exists a $G$-invariant symmetric} (resp. alternating) {\it nondegenerate $K$-bilinear form on $V$. Then there exists a $k$-bilinear form
on $V_k$ with the same properties.}

\smallskip 
 
 \smallskip
  \n  Assume now that $2$ is invertible in $R$, i.e. char$(k) \neq 2$, and hence char$(K) \neq 2.$ We may identify symmetric bilinear forms with quadratic forms. The symmetric part of th.A can then be restated, and made more precise, as follows:
      
   \smallskip
   {\bf Theorem B}. {\it Let $q$ be a nondegenerate $G$-invariant quadratic form on $V$. There exists a nondegenerate
   $G$-invariant quadratic form on $V_k$, whose class in the Witt ring $W(k)$ of  $k$  is the sum of the two Springer residues $\partial_1(q), \partial_2(q)$ of~$q$.} 
   
      \smallskip
   \n[For the definition of the Springer residues \ 
  $\partial_1, \partial_2 : W(K) \to W(k),$ see §3.3.]
   
   \smallskip
   
 As mentioned above, we shall prove a generalized form of these two theorems, the generalization consisting in replacing the action of a group by the action of an $R$-algebra with involution (th.5.1.4 and th.5.1.7).
 
 \smallskip
 
   The first three sections are about lattices; we define and give the main properties of what we call the {\it lower middle} and the {\it upper middle} of two lattices. The last two sections give the proof of theorems A and  B. 
   
   \smallskip
\n {\it Acknowledgment}.   I want  to thank R. Guralnick for his comments on a preliminary version of this paper, and for pointing out a paper of J. Thompson ([Th 84]), which contains all the necessary ingredients for proving th.A. We shall use some of Thompson's proofs in §4 and §5.
   
   \bigskip

 {\bf §0. Notation}

\medskip
\n 0.1. {\bf Lower middle and upper middle of a pair of integers}

\medskip

Let  $x,y$ be two integers.

\n The {\it lower middle} $m_-(x,y)$ of $(x,y)$ is defined as the largest integer $\leqslant \frac{x+y}{2}$ :
 
 \smallskip
   $ (0.1.1) \hspace{5mm}  m_-(x,y) = \lfloor \frac{x+y}{2} \rfloor.$
   
   \smallskip
\n    The {\it upper middle}  $m_+(x,y)$ of $(x,y)$ is the smallest integer $\geqslant \frac{x+y}{2}$ :
   
   \smallskip
   $(0.1.2) \hspace{5mm}  m_+(x,y) = \lceil \frac{x+y}{2} \rceil.$
   
   \smallskip

\n {\small  [This is similar to calling Wednesday and Thursday the {\it middle days} of the week.]}
   
   \smallskip
   
 \n  We have:
 
 \smallskip

    $(0.1.3) \hspace{5mm}  m_-(x,y)= \frac{x+y}{2} = m_+(x,y)$ \ if \ $ x \equiv y $ (mod 2)
     
   \smallskip

    $(0.1.4)  \hspace{5mm}  m_-(x,y)<\frac{x+y}{2} < m_+(x,y)$   \ and \  $ m_+(x,y) = m_-(x,y) + 1$  if  $ x \equiv y+1 $ (mod 2).

   \smallskip

   $(0.1.5) \hspace{5mm}  m_+(-x,-y) = - m_-(x,y)$.
   
   \smallskip

   $(0.1.6) \hspace{5mm}   m_-(x,y) = \sup_{n \in \Z} \ \inf(x-n,y+n)$.
   
   \smallskip

   $(0.1.7) \hspace{5mm}   m_+(x,y) = \inf_{n \in \Z} \  \sup(x-n,y+n)$.

\smallskip

$(0.1.8) \hspace{5mm}   m_-(x+1,y) = m_+(x,y)$ \ and \ $m_+(x+1,y) = m_-(x,y)+1$.

\medskip

   \n 0.2. {\bf Discrete valuations}

\medskip
We keep the same notation as in the introduction :  $K$  is a field with a discrete valuation $v : K^\times \to \Z$, which is extended to $K$ by putting $v(0) = + \infty$. The valuation ring $R$ is the set of all $x \in K$ with  $v(x) \geqslant 0$; the maximal ideal $\m$ of $R$ is the set of all $x \in K$ with $v(x) \geqslant 1$; we choose a generator $\pi$ of $\m$; we have
$v(\pi)=1$. The residue field is  $k = R/\m = R/\pi R$. 

The letter $V$ denotes a finite dimensional $K$-vector space.

\medskip

 {\bf §1. The lower and upper middles of a pair of lattices}

\medskip

\n 1.1. {\bf Definitions}

\smallskip

  Recall that a {\it lattice} of $V$ is a free $R$-submodule $L$ of $V$ such that the natural map   $K \otimes_R L \to V$ is an isomorphism. 
  
  Let $L$ and $M$ be two lattices. Then  $L \cap M$ is a lattice, and so is:
  
  \smallskip
  \hspace{15mm} $L+M =$ set of all $x+y$, with $x\in L, y\in M$. 
  
  \smallskip 
  We now define two lattices $m_-(L,M)$ and $m_+(L,M)$, which are sandwiched between  $L \cap M$ and
  $L+M$; we call them the  {\it lower middle} and {\it upper middle}  of  $L$ and $M$. (We shall see in prop.1.1.5 below how they are related 
  to the ``middles" of §0.1.) They are defined as follows :

   \smallskip
   $ (1.1.1) \hspace{5mm}  m_-(L,M) = \ $\mbox{\large{+}}$_{n\in \Z} \ (\pi^nL \   \cap \  \pi^{-n}M)$

\smallskip

$ (1.1.2) \hspace{5mm}  m_+(L,M) = \bigcap_{n\in \Z}\  (\pi^nL \ + \  \pi^{-n}M)$. 

\smallskip
\n  {\small[In these formulas, the \ $+$ \ symbol means `` submodule generated by ". For instance $m_-(L,M) $ is the $R$-submodule of $V$ generated by the lattices $ \pi^nL \   \cap \  \pi^{-n}M$. ]}
  
  \medskip
  
  Note that $L$ and $M$ play a symmetric role: we have  $m_{±}(L,M) = m_{±}(M,L)$.
  
  \medskip 
  
\n  {\it Remark.}\\ In (1.1.1) and (1.1.2), it is not necessary to run $n$ through the full set $\Z$; a finite subset suffices. For instance, in the case of (1.1.1), if  $n$  is large enough, then  $\pi^nL$  is contained in $L\cap M$,
 hence the $n$th term  $\pi^nL \   \cap \  \pi^{-n}M$  is contained in the $0$th-one; we may delete it from (1.1.1) without changing the sum. The same is true if $-n$ is large enough, by a similar argument. Same thing for (1.1.2). Hence the above infinite sums and intersections can be replaced
 by finite ones; this shows that  $m_+(L,M)$ and  $m_+(L,M)$ are lattices.
 
 As an example, suppose that $\pi^2L \subset M$ and $\pi^2M \subset L$. Then the terms with $|n| > 1$ can be deleted, and the formulas reduce to:
 
 \smallskip
   $m_-(L,M) = \ \ \ \pi L \cap \pi^{-1}M \ \ + \ \ L \cap M \ \  + \ \ \pi^{-1}L \cap \pi M$,
   
   \smallskip
   
   $m_+(L,M) = \ (\pi L + \pi^{-1}M) \ \cap  (L + M) \ \cap \ (\pi^{-1}L + \pi M)$. 
   
   \smallskip
   More generally, if $a \geqslant 0$ is such that  $\pi^aL \subset M$ and $\pi^aM \subset L$, the terms with $|n| > a/2$ can be deleted.
   
   \medskip
   
   \n {\it Example $:$ middles of twisted lattices}
   
   \smallskip
   
     If $L$ is a lattice, and if $a \in \Z$, the $a$-{\it twist} $L(a)$ of $L$ is defined as :
     
      \smallskip
     $ (1.1.3)  \hspace{5mm} L(a) = \pi^{-a}L $.
     
      \smallskip
   
   We have  $a \leqslant b \Rightarrow L(a) \subset L(b)$, and:
   
    \smallskip
   $ (1.1.4)  \hspace{5mm}  L(a) + L(b) = L(\sup(a,b)) $ \ and \ $L(a) \cap L(b) = L(\inf(a,b))$.
   
   \medskip
   
 \n  {\bf Proposition 1.1.5}. {\it If $x,y \in \Z$, then} :
   
    \smallskip
   $(1.1.6) \hspace{5mm} m_-(L(x),L(y)) = L(m_-(x,y)),$
   
    \smallskip
   $(1.1.7) \hspace{5mm} m_+(L(x),L(y)) = L(m_+(x,y))$
   
    \smallskip
   {\small [In other words : on the set of all twists of a given lattice, the ``middle" operations
   coincide with those defined on $\Z$ in §0.1.]} 
   
    \smallskip
    
    \n {\it Proof.} \\ This follows from (1.1.4), combined with (0.1.6) and (0.1.7).

   \bigskip

\n 1.2. {\bf Basic properties of the lower and upper middles of two lattices}

\medskip

  As above, let $L$ and $M$ be two lattices of $V$.
  
  \medskip
   
  \n   {\bf Proposition 1.2.1}. {\it We have }:
     
     \smallskip
     
  $   (1.2.2) \hspace{5mm} L \cap M \ \subset \ m_-(L,M) \ \subset m_+(L,M) \ \subset \ L+M.$
   
   \smallskip
   
  \n  {\it Proof.}\\
      The inclusions $L \cap M \ \subset \ m_-(L,M)$  and   $m_+(L,M) \ \subset \ L+M$ are clear. 
   
 \n  Let us show that
    $m_-(L,M) \ \subset m_+(L,M)$. Note first that we have :
    
    \smallskip
    
      $(1.2.3)  \hspace{5mm} \pi^ aL \cap \pi^{-a}M \ \subset \  \pi^bL + \pi^{-b}M$ \quad for every $a,b \in \Z$. 
    
    \smallskip
 \n   Indeed, if $a \geqslant b$, this follows from :
      
      \smallskip
     \hspace{5mm}  $\pi^ aL \cap \pi^{-a}M  \ \subset \pi^aL \ \subset \pi^bL \ \subset \ \pi^bL + \pi^{-b}M$;
      
      \smallskip
      
   \n   similarly, if $b \geqslant a$, (1.2.3) follows from :

       \smallskip
  \hspace{5mm}     $\pi^ aL \cap \pi^{-a}M  \ \subset \pi^{-a}M  \ \subset \pi^{-b}M \ \subset \ \pi^bL + \pi^{-b}M$.
      
      \smallskip
      
  \n    Since $m_-(L,M)$ is generated by the $\pi^ aL \cap \pi^{-a}M$, formula (1.2.3) shows that it is contained in every 
   $\pi^bL + \pi^{-b}M$, hence also in their intersection, which is $m_+(L,M)$.
   
   \medskip
   
   \n {\it Remark}. \\ The $R$-modules $(L+M)/m_+(L,M)$ and $m_-(L,M)/(L \cap M)$ are (non canonically) isomorphic; this is proved by the method of reduction to dimension~1 used in the proof of th. 1.2.4 below. \\
   Similarly, we have $(L+M)/m_-(L,M) \ \simeq \ m_+ (L,M)/(L \cap M)$.
   
   \medskip
    \n   {\bf Theorem 1.2.4}. {\it We have}  $\pi.m_+(L,M) \ \subset \ m_-(L,M)$.
   
   \smallskip

   \n {\it Proof.}\\ Suppose first that  $\dim V = 1$. In that case, all the lattices of $V$ are twists of one of them. By prop.1.1.5, the formula  $\pi.m_+(L,M) \ \subset \ m_-(L,M)$ is equivalent to the obvious formula \ $-1+m_+(x,y) \leqslant m_-(x,y)$ \ for $x,y \in \Z$.
     
     We now reduce the general case to that one. To do so, let $V = \oplus \ V_i$ be a splitting of $V$ as a direct sum of subspaces $V_i$. We say that this splitting is {\it compatible} with a lattice $P$  if $P = \oplus \ P_i$, where $P_i = P \cap V_i$. If $V = \oplus \ V_i$ is compatible with two lattices  
     $L$ and $M$, then the same is true for $L \cap M$, and we have $(L \cap M)_i = L_i \cap M_i$; same for $L + M$, $\pi^ aL \cap \pi^bM$, $m_±(L,M)$. Hence, {\it if th. 1.2.4 is true for the $V_i$, it is true for $V$}. 
   
It only remains to prove :
    
      \smallskip
     
     \n   {\bf Lemma 1.2.5}. {\it There exists a splitting $V = \oplus \ V_i$, with $\dim V_i =1$ for every $i$, which is compatible with both $L$ and $M$.}
     
     \smallskip
   \n {\it Proof.}\\ By replacing $M$ by a suitable $\pi^nM$, we may assume that $M \subset L$. Since $R$ is a principal ideal domain, there exists an $R$-basis $(x_i)$ of $L$, and nonzero elements $a_i$ of $R$ such that the $(a_ix_i)$ make up a basis of $M$, cf. e.g. [A VII, §4, th.1]. The splitting $V = \oplus \ Ke_i$ has the required properties: one has  $L = \oplus \ Re_i$ and $M = \oplus \ Ra_ie_i$.

  \medskip
  
  \n {\it Remark. }\\ Let us mention three formulas which can also be proved by reduction to dimension 1:
  
    \smallskip
  (1.2.6) \hspace{5mm} $m_-(\pi L,M) = \pi.m_+(L,M)$,
  
    \smallskip
  (1.2.7) \hspace{5mm} $m_+(\pi L,M) = m_-(L,M),$
  
  \smallskip
  
  (1.2.8) \hspace{5mm} $m_±(L,M) = m_±(L+M, L\cap M)$.
    
  \medskip 
  \n {\it Building interpretation.}
  
  \smallskip
  
   When $R$ is complete for the $\pi$-adic topology, the lattices of  $V$  may be viewed as the vertices of the {\it affine building} $X$ associated by Goldman-Iwahori to $\GL(V)$, cf. [GI 63] and [Ge 81]; the space $X$ is isomorphic to the product of {\bf R} by the Bruhat-Tits building of $\SL(V)$, cf. [BT 84]; the apartments of $X$ correspond to the splittings of  $V$ as direct sums of 1-dimensional subspaces, and lemma 1.2.5 is a special case of the fact that any two points of $X$ are contained in an apartment, cf. [Ge, 2.3.4]. If $L,M$ are two lattices, and $[L],[M]$ are the corresponding  vertices of  $X$, any barycenter $x[L] + y[M]$, with $x,y$ real $\geqslant 0$, $x+y=1$, makes sense as a point of $X$ (such barycenters make up the {\it geodesic segment} joining $[L]$ to $[M]$). In particular,
  $\frac{1}{2}$$[L]$+ $\frac{1}{2}$$[M]$ makes sense; it is the middle (in the standard meaning of the word) of that geodesic segment. This middle point is not always a vertex, i.e., it does not always correspond to a lattice\footnote{But it does after extension of scalars to $K(\sqrt{\pi})$.}. One can check that it is indeed a vertex if and only if $m_+(L,M)$ and $m_-(L,M)$ coincide, and in that case we have the good-looking formula:  
  
  \smallskip
  \hspace{2cm} $\frac{1}{2}$$[L]$+ $\frac{1}{2}$$[M]$ = $[m_+(L,M)] = [m_-(L,M)] .$

   \bigskip

\n 1.3. {\bf The lower and upper middle submodules of a torsion $R$-module}

\smallskip

\n {\small [The content of this section will not be used in the rest of this paper.]}

\medskip 
Let $T$ be a torsion $R$-module of finite exponent, i.e., an $R$-module such that there exists  $N \geqslant 1$ with
$\pi^NT = 0$. Let us define its {\it lower middle} and {\it upper middle} submodules by the following formulas :

  \smallskip
(1.3.1) \hspace{5mm}  $m_-(T) = \ $\mbox{\large{+}}$_{n\geqslant 0} \ (\I \ \pi_T^n \ \cap \  \Ker \pi^n_T)$,

  \smallskip
(1.3.2) \hspace{5mm}  $m_+(T) = \bigcap_{n\geqslant 0} \ (\I \ \pi_T^n \ + \  \Ker \pi^n_T)$,

  \smallskip
  
  \n where  $\pi_T$ is the endomorphism of $T$ defined by  $\pi$. In these formulas, it is enough to let $n$ run from $1$ to $N-1$:
  larger $n$'s contribute nothing, since then $\pi^n_T$ = 0, and $n=0$ does not give anything either.
  
  \smallskip
  
  The properties of $m_-(T)$ and $m_+(T)$ are analogous to the ones of §1.1 and §1.2. One may sum them up by :
  
  \medskip
  
     \n   {\bf Proposition 1.3.3}. {\it We have}  :
     
     \smallskip
    (1.3.4) \hspace{5mm} $\pi.m_+(T) \subset m_-(T)  \subset  m_+(T)$.
    
    \smallskip
    (1.3.5) \hspace{5mm} $m_+(T) \ \simeq  \ T/m_-(T).$ 
         
    \smallskip
    (1.3.6) \hspace{5mm} $m_-(T) \ \simeq \ T/m_+(T).$ 
    
    \medskip
    
    The proof is analogous to that of th.1.2.4: check first the case where $T$ is cyclic,
    i.e., $T \simeq R/\pi^mR$ for some $m \leqslant N$, and then use the well known fact that $T$ is a direct sum (finite or infinite) of cyclic submodules.
    
    \medskip
    
    This construction is related to that of §1.1 and §1.2: if $L$ and $M$ are two lattices of $V,$ we have  $m_±(T)= m_±(L,M)/(L \cap M)$, where  $T = (L+M)/(L\cap M)$.
    
    \smallskip
    
    \n{\small {\it Application to isogenies between abelian varieties}
    
    Let  $A$ and $B$ be abelian varieties over a field  $F$, and let  $\varphi : A \to B$ be an isogeny. Assume that deg$(\varphi)$ is a power of a prime number $\ell$ distinct from char$(F)$. Let $F^s$ be a separable closure of $F$, and let $T = \Ker \ A(F^s) \to B(F^s)$. By assumption, $T$ is a finite abelian $\ell$-group of order deg$(\varphi)$, hence is a torsion module over $\Z_\ell$. Let $m_-(T), m_+(T)$ be the middles of  $T$. The Galois group  $\Gal(F^s/F)$ acts on $T$, and stabilizes $m_-(T)$ and  $m_+(T)$. Thus, there is a factorization of $\varphi$ as  $A \to  A_-  \to A_+ \to B$, where  $A_-$ and $A_+$ are abelian varieties over $F$, and the kernel of $A(F^s) \to A_±(F^s)$ is $m_±(T)$; by $(1.3.4)$, the kernel of $A_-  \to A_+$ is killed by $\ell$.
    
    An interesting special case is when $B$ is the dual variety of $A$, and $\varphi$ is a polarization of degree a power of $\ell$. In that case, one can show that  $A_-$ and $A_+$ are dual of each other. Hence, {\it if an abelian variety over $F$
    has a polarization of $\ell$-power degree, there is an $F$-isogenous one which has the same property, and for which the polarization kernel is killed by}  $\ell$.}

     \bigskip
    {\bf §2. The Brauer-Nesbitt theorem}
    
    \smallskip
    
  \n2.1. {\bf Notation}
    
       Let $A$ be an $R$-algebra and let  $(a,x) \mapsto ax$ be an $R$-bilinear map  $A \times V \to V$ which makes $V$ into an $A$-module. 
       
      \n   {\small[If $A$ is the group algebra $R[G]$ of a group $G$, this means that $G$ acts linearly on $V$.]}
       
        \smallskip

       \n {\bf Lemma 2.1.1}. {\it The following properties are equivalent}:
       
       \smallskip 
       (2.1.2) {\it There exists a lattice $L$ of $V$ with is stable under the action of  $A$, i.e.,} $A.L = L$.
       
       (2.1.3) {\it The image of $A$ in $\End(V)$ is a finitely generated $R$-module} .
      
      \smallskip
      
      \n {\it Proof.} 
      
       Let $A_V$ be the image of $A$ in $\End(V)$. If $L$ is as in (2.1.2), then we have $A_V \subset \End_R(L)$, and (2.1.3) follows. Conversely, if $A_V$ is finitely generated, let $M$ be a lattice, and let $L = A_V.M$; then $L$ is a finitely generated $R$-module which contains $M$; hence it is a lattice such that $A.L = L$. 
      
      \smallskip
      
      If (2.1.2) and (2.1.3) hold, we shall say that the action of $A$ on $V$ is {\it bounded}. A lattice $L$ with $A.L =L$ will be called an {\it $A$-lattice}.
      
      \medskip   
      
    \n  2.2. {\bf Semisimplification}
    
    \smallskip
      Let $L$ be an $A$-lattice; then $E_L = L/\pi L$ is a module over the $k$-algebra $A_k = A/\pi A$. Let $S_1,...,S_m$ be the successive quotients of a Jordan-Hölder filtration of $E_L$; the direct sum  $ E_L^{\s} = S_1   
  \oplus ... \oplus S_m$ is a semisimple
      $A_k$-module which is independent (up to isomorphism) of the chosen Jordan-Hölder filtration; that module is called the {\it semisimplification} of  $E_L$.       
      
      \medskip
      
  \n {\bf Theorem 2.2.1} (Brauer-Nesbitt). {\it Let $L$ and $M$ be two $A$-lattices of $V$. Then the  $A_k$-modules    $ E_L^{\s} $
      and $ E_M^{\s}$  are isomorphic.} 
      
      \smallskip
      This was proved by Brauer and Nesbitt when $A$ is the $R$-algebra of a finite group. The proof in the general case is similar. There are two steps:
      \smallskip
      
    \n  (2.2.2) {\it The special case where $\pi L \subset M \subset L$.}  
      
      \medskip
  Let $T = L/M$. We have an exact sequence of $A_k$-modules :
  
  \smallskip
  (2.2.3)  \hspace{15mm} $0 \to T \to E_M \to E_L \to T \to 0$,
  
  \smallskip 
\n where the map  $T \to E_M$ is induced by  $x \mapsto \pi x$. This implies an isomorphism :

\smallskip
 (2.2.4)  \hspace{15mm} $ T^{\s} \oplus E_L^{\s}  \simeq  E_M^{\s} \oplus T^{\s}$;
 
 \smallskip
\n hence   $ E_L^{\s}  \simeq  E_M^{\s} $. 

\medskip

\n (2.2.5) {\it Reduction of the general case to the special case}.

\smallskip
   Replacing $M$ by a multiple  $\pi^aM$ does not change $E_M^{\s}$; hence we may assume that 
   $M \subset L$. There exists $n \geqslant 0$ such that $\pi^nL \subset M$; choose such an $n$ and use induction on $n$. The case $n=0$ is trivial since $M=L$; if $n > 0$, define $N = \pi^{n-1}L + M$. Since $\pi^{n-1}L \subset N \subset L$, the induction hypothesis implies that
$ E_L^{\s}  \simeq  E_N^{\s} $. Since $\pi N \subset M \subset N$, we have $ E_M^{\s}  \simeq  E_N^{\s} $ by (2.2.2). Hence  $ E_L^{\s}  \simeq  E_M^{\s} $.

\medskip

\n {\it Alternate proof}.\\ Since $E_L^{\s} $
      and $ E_M^{\s}$ are semisimple, to prove that they are isomorphic it is enough to show that, for every $a \in A$, the characteristic polynomials of  $a$, acting on these two modules, are the same (this criterion, for group algebras of finite groups, is due to Brauer - for the general case, see Bourbaki [A VIII.377, §20, th.2]); but this is clear, since these polynomials
      are the reduction mod $\pi$ of the characteristic polynomial of $a$ acting on $V$.

\smallskip
\n {\it Notation.} The $A_k$-module $E_L^{\s} $ will be called the {\it reduction} mod $\pi$ of $V$; we shall denote it by $V_k$.
It is defined up to a non canonical isomorphism.

\bigskip

 {\bf §3. Bilinear forms and lattices}

\smallskip

Let $B(x,y)$ be a nondegenerate $K$-bilinear form on $V$, which is $\e$-symmetric, with $\e= ±1$, i.e., $B(x,y)= \e B(y,x)$ if $x,y \in V$.

\medskip

\n 3.1. {\bf Dual lattices}

\smallskip

   Let $L$ be a lattice of $V$, and let $L'$ be the {\it dual lattice} of $L$, namely the set of all $x\in V$
  such that $B(x,y) \in R$
  for every  $y \in L$. If $x \in L'$, the map $ y \mapsto B(x,y)$ is $R$-linear; this gives an isomorphism
  $L' \to \Hom_R(L,R)$; we may thus identify $L'$ with the usual ``$R$-dual'' of $L$.  We have $(L')' = L$.
  \smallskip
  
  A lattice  $L$  is {\it self-dual} (or {\it unimodular}) if $L'=L$; if $(e_1,...,e_n)$ is an $R$-basis of $L$, this means that both the matrix $ \sf{B}$ $=(B(e_i,e_j))$ and its inverse have coefficients in $R$, i.e., we have  $ \sf{B}$ $ \in {\rm GL}_n(R).$ 
  
  \smallskip
    We say that $L$ is {\it almost self-dual} if $\pi L' \subset L \subset L'$, i.e., if the matrices $\sf{B}$  et $\pi\sf{B}^{-1}$ have coefficients in $R$.
    
    \medskip
    
     \n {\bf Theorem 3.1.1}  {\it Let $L$ be a lattice and let $L'$ be its dual. Then the lower middle $m_-(L,L')$ of $L$ and $L'$ is an almost self-dual lattice whose dual is} $m_+(L,L')$.
     
     \smallskip 
     \n {\it Proof}.\\ The functor ``lattice $\mapsto$ dual lattice'' transforms finite intersections into finite sums, and conversely. By (1.1.1) and (1.1.2), this shows that  $m_-(L,L')$  and  $m_+(L,L')$ are dual of each other. By prop.1.2.1 and th.1.2.4, we have $$\pi.m_+(L,L') \subset m_-(L,L') \subset m_+(L,L'),$$ 
     hence $m_-(L,L')$ is almost self-dual.
    
    \medskip 
    
\n 3.2. {\bf The bilinear forms $b_1$ and $b_2$ associated with an almost self-dual lattice}    

Let  $L$  be an almost self-dual lattice. The inclusions $\pi L' \subset L \subset L'$ give an exact sequence of $k$-vector spaces :

 \smallskip
(3.2.1) \quad $ 0 \to L'/L \to L/\pi L \to L/\pi L' \to 0, $ 

 \smallskip
  \n  where the map  $L'/L \to L/\pi L$ is given by  $x \mapsto \pi x$.
    
     \smallskip
     
    By passage to quotients,  the bilinear form $B$ gives a $k$-bilinear forms  $b_1$  on $ L/\pi L'$; similarly, $\pi B$ defines a $k$-bilinear form $b_2$ on $L'/L$. These forms are $\e$-symmetric, and nondegenerate.

     \medskip

    \n 3.3. {\bf The quadratic case, and the Springer residues} 
    
    Let us assume now that $\e = 1$, i.e., that $B$ is symmetric. Assume also that $\car(k) ≠ 2$, i.e., that 2 is invertible in $R$. We may thus identify symmetric bilinear forms and quadratic forms. Let $q(x) = B(x,x)$ be  the quadratic form
    defined by $B$, and let  $[q]$ be its image in the Witt ring  $W(K)$ of the field  $K$ (for the definition and basic properties of the Witt ring, see [La 05, chap.I-II]). Similarly, let  $q_1, q_2$ be the quadratic forms defined by the $k$-bilinear forms  $b_1,b_2$ of §3.2, and let $[q_1], [q_2]$ be their images in the Witt ring $W(k)$. 
    
      \smallskip
      Recall (cf. [Sp 55], [La 05, chap.VI]) that Springer has defined two ``residue'' maps
      
   \smallskip   
 \hspace{2cm}     $\partial_1, \partial_2 : W(K) \ \to \ W(k).$
      
      \smallskip
      
          The map  $\partial_1$ is a ring homomorphism; the map $\partial_2$ is additive (and depends on the choice of the uniformizing element  $\pi$). They are characterized by these properties, together with their values for 1-dimensional
    quadratic forms, which are as follows :
    
    \smallskip 
    
  (3.3.1)   If  $u \in R^\times$ has image  $\overline{u}$ in $ k$, then the images  of the $1$-dimensional form
    $\langle u \rangle$ by $\partial_1,  \partial_2$ are $\langle \overline{u} \rangle,  0$, and those of $\langle u\pi \rangle$
    are $0,  \langle \overline{u} \rangle$.
    
    \smallskip
    The map \  $(\partial_1, \partial_2) : W(K) \to W(k) \times W(k)$ \ is surjective; it is bijective if $K$ is complete.

    \medskip
    
     \n {\bf Theorem 3.3.2.}  {\it Let $L$ be an almost self-dual lattice, and let $q_1,q_2$ be the corresponding quadratic forms over $k$. Then $\partial_1([q]) = [q_1]$ and $\partial_2([q]) = [q_2]$.}

      \medskip
      
      \n {\it Proof.}\\ We use the same method as for th.1.2.4, namely reduction to dimension~1 (in which case the formulas are obvious). What is needed is the following orthogonal analogue of lemma 1.2.5 :
      
      \smallskip
      
       \n {\bf Lemma 3.3.3.}  {\it If $M$ is any lattice of $V$, there exists an orthogonal splitting \ $V = \oplus \ V_i$, with $\dim V_i = 1$ for every $i$, which is compatible with $M$.}
       
       \n [The orthogonality assumption means that $V_i$ and $V_j$ are orthogonal
       for the bilinear form $B$ if  $i ≠ j.$ It implies that the splitting is compatible with the dual $M'$ of $M$. ]
       
    \smallskip  
    
      \n {\it Proof of the lemma.} 
      
      Use induction on $\dim V$. Let $m = \inf_{x \in M} v(q(x))$; we have $m > - \infty$ . Choose  $x \in M$ with  $v(q(x))=m$.
      For $y,z \in M$, we have  $v(B(y,z)) \geqslant m$ : this follows from the formula  $2B(y,z) = q(y+z)-q(y)-q(z)$ since the valuations of  $q(y+z), q(y), q(z)$ are $\geqslant m$. If $y \in V$, put $\ell(y) = B(x,y)/B(x,x)$. The linear form
      $\ell : V \to K$ is such that  $
     \ell(x) = 1$, and it maps $M$ onto $R$. It thus gives a splitting of $V$ as $ Kx \oplus \Ker \ell$, namely $y \mapsto (\ell (y)x, y -\ell(y)x$; this splitting is compatible with $M$. It is an orthogonal splitting, since  $B(y - \ell (y)x,x) = 0$ for every  $y$. The lemma follows by applying the induction assumption to the vector space $\Ker \ell$ and its lattice $ M \cap \Ker \ell$.

   \medskip
    
      \n {\bf Corollary 3.3.4.}  {\it The classes in $W(k)$ of the two quadratic forms associated with an almost self-dual lattice do not depend on the choice of that lattice.}
    
 \medskip
      \n {\it Remark}. The fact that the class of $q_1$ in  $W(k)$  is the same for any two almost self-dual lattices $L$ and $M$ does not imply that  these two quadratic forms have the same dimension: they may differ by hyperbolic factors. But, if they do have the same dimension, a theorem of Bayer-Fluckiger and First shows that $L$ and $M$ are {\it isomorphic as quadratic $R$-modules} ([BF 17], th.4.1); hence they are conjugate of each other by an element  of the orthogonal group of $(V,q)$. 
      
\bigskip

 {\bf §4. Semisimplification of symplectic and quadratic modules over an algebra with involution}

\smallskip
This section is essentially independent of §§1,2,3: we work over a field, and not over a discrete valuation ring.

\medskip

 \n 4.1. {\bf Notation}

\smallskip
  Let  $k$ be a field, and  let $A_k$ be a $k$-algebra with a $k$-linear involution denoted by  $a \mapsto a^*$; we have $a^{**} = a$ and $(ab)^* = b^*a^*$ for every  $a,b \in  A_k$. 
  
  \smallskip
   Let  $E$ be an $ A_k$-module which is finite dimensional over $k$ (i.e., a ``linear representation'' of  $ A_k$), and let $b$
   be a nondegenerate bilinear form on $E$ with the following property:
   
   \smallskip 
   (4.1.1) $b(ax,y) = b(x,a^*y)$ for every  $a \in  A_k, \  x,y \in E.$
   
   \smallskip

   We then say that $b$ is {\it compatible} with the $ A_k$-module structure of $E$. When $ A_k$ is the group algebra $k[G]$ of a group $G$, with its canonical involution  $g^*=g^{-1}$ for every $g\in G$, this means that $b$ is invariant by  $G$.
   
    \smallskip
   We assume one of the following :
   
    \smallskip
   (4.1.2)  $b$ is alternating.
   
    \smallskip
   
   (4.1.3) $b$ is symmetric and $\car(k) ≠ 2$.

    \smallskip
   In the first case, $E$ is called a {\it symplectic $ A_k$-module}, and in the second case, it is called an {\it orthogonal $ A_k$-module}.

\medskip

\n 4.2. {\bf Semisimplification}

\smallskip

Let $ A_k, E, b$ be as above, and let $E^{\s}$ be the semisimplification of the $ A_k$-module $E$.

\medskip

     \n {\bf Theorem 4.2.1.}  {\it In case $(4.1.3)$ $($resp. in case $(4.1.2))$, there exists a symmetric $($resp. alternating$)$
     bilinear form on  $E^{\s}$ with the following two properties}:
     
     \smallskip
     
     (4.2.2) {\it It is compatible with the $ A_k$-module structure of  $E^{\s}$}.
     
     \smallskip
     
     (4.2.3) {\it It is isomorphic to $b$.}

\smallskip

\n One may sum up (4.2.2) by saying that the semisimplification of a quadratic (resp. alternating) module is quadratic (resp. alternating). As for (4.2.3), it is obvious in the symplectic case, since all nondegenerate 
alternating forms of a given rank are isomorphic; in the orthogonal case, it means that the corresponding quadratic forms have the same class in the Witt ring  $W(k)$.

\smallskip

\n {\it Proof.}\\ We use the method of [Th 84, §2]. Let $S$ be an $A_k$-submodule of $E$, which is totally isotropic  for $b$ (i.e. $b(x,y)=0$ for all $x,y \in E$), and is maximal for that property. Let $S_\bot$  be its orthogonal relative to $b$; because of (4.2.2), it is a submodule of $E$, and we have $ 0 \subset S \subset S_\bot \subset E$. The form $b$ defines a nondegenerate form $b_1$ on $S_\bot/S$.

\smallskip

   \n {\bf Lemma 4.2.4.} (i) {\it The $A_k$-module $X= S_\bot/S$ is semisimple, and the only totally isotropic submodule of $X$ is $0$.}
   
   (ii) {\it In the orthogonal case $(4.1.3)$, the form  $b$ is isomorphic to the direct sum of  $b_1$ and an hyperbolic
   form of rank $2 \dim S$.}

\smallskip
\n {\it Proof of} (i).\\ If $Y$ is a totally isotropic submodule of $X$, its inverse image in $S_\bot$ is totally isotropic, hence equal to  $S$, i.e., $Y$= 0. If $Z$ is a submodule of $X$, then $Y = Z \cap Z_\bot$ is totally isotropic, hence $0$; we have $X = Z \oplus Z_\bot$; this shows that every submodule of $X$ is a direct summand, i.e., $X$ is semisimple.

\smallskip

\n {\it Proof of} (ii).\\ Let $M$ be a subvector space of  $S_\bot$ such that $S_\bot =S \oplus M$. The restriction of $b$ to  $M$
   is nondegenerate, and the projection $M \to S_\bot/S$ is an isomorphism of quadratic spaces; we have $E = M \oplus M_\bot$, and  $\dim M_\bot = 2 \dim S.$ The form $b$ splits as $ b_1 \oplus h$, where $h$ is the quadratic form of $M_\bot$; the form $h$ is hyperbolic, since that space contains the totally isotropic subspace  $S$, of dimension  $ \frac{1}{2}$ $\dim M$. This proves (ii).

   \smallskip
   
   \n {\it End of the proof of th.4.2.1}.\\ The bilinear form  $b$  defines a duality between  $S$  and  $E/S_{\bot}$;
   hence we may identify $E/S_{\bot}$ with the $k$-dual $S'$ of $S$ (with its natural $A_k$-structure). We have
   $$ E^{\s} = X \oplus (S^{\s} \oplus S'^{\s}),$$
   \n since $X$ is semisimple, cf. 4.2.4 (i). It is easy to see that $S'^{\s}$ is isomorphic to the dual of $S^{\s}$.
   Hence, $ E^{\s}  \simeq X \oplus (Y \oplus Y')$, where  $Y = S^{\s}$. We then put on $ E^{\s} $ the bilinear
   form which is the direct sum of the form  $b_1$ on $X$ and the natural bilinear form (symmetric or alternating, as needed) on $Y \oplus Y'$. By lemma 4.2.4, that form has properties (4.2.2) and (4.2.3).

   \bigskip

 {\bf  §5. Proof of theorem B}

\smallskip

 We now prove the theorems stated in the introduction; the proofs will merely consist in putting together the results of §§2,3,4.

\medskip
\n 5.1. {\bf The setting}

\smallskip
  We go back to the standard notation $(K,R,V)$ and we assume that  $V$ is a module over an $R$-algebra $A$ with involution. We also assume that this action is {\it bounded} (cf. §2.1), i.e., that there exists a lattice of $V$  which is
  stable under $A$.
    
   Let $B$ be a nondegenerate $K$-bilinear form on $V$, which is compatible with the action of $A$; this means (as in (4.1.1)) :
  
  \smallskip
  
  (5.1.1)  $ B(ax,y) = B(x,a^*y)$  \  for every  \  $a\in A, \ x,y \in V.$
  
  \smallskip
  As in §4, we assume one of the following:
  
   \smallskip
   
   (5.1.2)  $B$  is alternating.
   
    \smallskip
    
    (5.1.3) $B$ is symmetric, and  $\car(k) ≠ 2$.
  
   \bigskip
   
   By the Brauer-Nesbitt theorem (see §2.2), the semisimplification $V_k$ of  $V$  is well-defined; it is
   a module over the $k$-algebra $A_k = A/\pi A$.
   
   \medskip

   \n {\bf Theorem 5.1.4.}  {\it There exists a nondegenerate $k$-bilinear form $b$ on $V_k$ with the following properties}:
   
   \smallskip
   
   (5.1.5) {\it It is compatible with the action of $A_k$, i.e., it satisfies condition} (4.1.1).
   
   \smallskip 
   
   (5.1.6)  {\it It is alternating $($resp. symmetric$)$ \it if $B$ is.}
   
   \smallskip
   
  \n {\small  [More shortly : if $V$ is a symplectic (resp. orthogonal) module, so is $V_k$.]}
  
  \smallskip
  
   Note that condition (5.1.5) alone would be easy to satisfy: since $V$ is isomorphic to its dual, the same is true for $V_k$, and that is equivalent to the existence of a bilinear form $b$ compatible with the action of $A_k$.
      
   \medskip
   
   In the orthogonal case (5.1.3), one may ask more of the form  $b$. To state it, let us denote by $q$ (instead of $[q]$)
   the Witt class in $W(K)$ of the quadratic form  $q(x) = B(x,x)$, and let $q_1,q_2 \in W(k)$ be its images by the Springer residue maps, cf. §3.3. Then:
   
   \medskip
   
   \n {\bf Theorem 5.1.7.}  {\it The form $b$ of th.5.1.4 can be chosen to have the following property} :
   \smallskip
   
   (5.1.8) {\it There exists an $A_k$-orthogonal splitting $V_k = E_1 \oplus E_2$ such that the Witt class of the restriction to $E_i$ \ $(i=1,2)$ of $x \mapsto b(x,x)$ is } $q_i$.
   
   \medskip
   
   Note that (5.1.8) implies that the Witt class of the quadratic form  $b(x,x)$ on  $V_k$ is $q_1+q_2$. When
   $A$ is a group algebra $R[G]$, we recover th.B of the Introduction. Similarly, th.5.1.4 applied to $R[G]$ gives th.A.
    
      \medskip
      
      Theorems 5.1.4 and 5.1.7 will be proved in §5.3.
      
      \medskip
      
       \n   5.2. {\bf Existence of almost self-dual $A$-lattices}

      \medskip

        \n {\bf Theorem 5.2.1.}  {\it There exists an $A$-lattice of $V$ which is almost self-dual} (cf. §3.1) {\it relative to the bilinear form $B$.} 
   
   \smallskip
   
   \n {\it Proof.}\\ Let $L$ be an $A$-lattice of $V$. Condition (5.1.1) shows that the dual $L'$ of $L$ is also an $A$-lattice, and so are the lattices  $\pi^nL, \pi^nL'$, and hence also the lower middle $m_-(L,L')$, which is almost self-dual by th.3.1.1.
   
       \medskip
       
       \n {\it Alternate proof}.\\
        The proof above uses the ``middle'' notions of §§1,2. Here is a direct proof,
   taken from [Th 84]:

   Choose an $A$-lattice  $L$, with $L \subset L'$, and which is maximal for that property. Let  $m$  be the smallest
   integer such that $\pi^mL' \subset L$. If $m = 0$, or $1$, then $L$ is self-dual, or almost self-dual. Let us show
   that $m \geqslant 2$ is impossible. Define $M =\pi^{m-1}L' + L$; we have $M' = \pi^{1-m}L \cap L'$. The inequality     $m \geqslant 2$ implies $\pi^{m-1}L' \subset   \pi^{1-m}L$, hence $\pi^{m-1}L' \subset M'$; since $L \subset M'$, this shows that  $M$ is contained in $ M'$. The maximality of  $L$  then implies $M = L$, hence $\pi^{m-1}L' \subset L$, which contradicts
   the minimality of $m$.
   
   \smallskip
   
   \n {\it Remark}. This short proof is not in fact very different from the first one. Indeed, it amounts to construct an
   almost self-dual $A$-lattice, starting with any lattice $L$ contained in its dual, by choosing  $m$  with $\pi^mL' \subset L$, then replacing  $L$  by  $\pi^{m-1}L' + L$, and iterating until one gets $m = 0$ or $m=1$. But, if one writes
   down the end result of that process, a simple computation shows that one finds the lower middle lattice $m_-(L,L')$,
   exactly as in the first proof. 
   
     There is however a definite advantage in using a maximal lattice : in that case, $L'/L$ does not contain any nonzero
     totally isotropic submodule, hence it is a semi-simple $A_k$-module, cf. lemma 4.2.4. This simplifies the proofs  of the next section.

     \medskip
    \n   5.3. {\bf Proof of th.5.1.4 and th.5.1.7}
      
      \medskip
   Let $L$ be an almost self-dual $A$-lattice, cf. th.5.2.1. Let $F_1 =L/\pi L', \ F_2 = L'/L$.  By (3.2.1), we have an exact sequence of $A_k$-modules :
   
   \smallskip
   
   (5.3.1)   \  $0 \to F_2 \to L/\pi L \to F_1 \to  0. $ 
   
   \smallskip
   
   Let $E_1= F_1^{\s}, E_2 = F_2^{\s}$, and let $V_k = (L/\pi L)^{\s}$. The exact sequence above gives a splitting :
   
   \smallskip
   
   (5.3.2)  \ $V_k = E_1 \oplus E_2$.
   
   \smallskip 
   
   As explained in §3.2, the bilinear form $B$ defines $k$-bilinear forms $b_1$ and $b_2$ on $F_1$ and $F_2$; these forms are compatible with the action of $A_k$ and they are alternating (resp. symmetric) if $B$ is alternating (resp.
   symmetric). By th.4.2.1, applied to $F_1$ and $F_2$,  there exist $A_k$-compatible forms $b'_1, b'_2$ on $E_1,E_2$ which are isomorphic (as
   bilinear forms) to $b_1,b_2$; in the orthogonal case, th.3.3.1 shows that the Witt classes of these forms are the Springer residues of $q$. Using (5.3.2), we define a bilinear form on $V_k$ as the orthogonal sum of $b_1'$ on $E_1$ and $b_2'$ on $E_2$. All the properties of th.5.1.4 and th.5.1.7. hold.
   
   \smallskip
   
   This concludes the proof.

   \begin{center}
  {\bf References}   
   \end{center}

   \bigskip
   
   [A VII] N. Bourbaki, {\it Algèbre$,$ Chapitre $7,$ Modules sur les anneaux principaux}, Paris, Masson, 1981; English translation, {\it Algebra} II, Springer-Verlag, 1989.
   
   [A VIII] -------, {\it Algèbre$,$ Chapitre $8,$ Anneaux et modules semi-simples}, new revised edition, Springer-Verlag, 1998.
   
   [BF 17] E. Bayer-Fluckiger \& U.A. First, {\it Rationally isomorphic hermitian forms and torsors of some non-reductive groups}, Adv. Math. 312 (2017), 150-184.
   
  [BT 84] F. Bruhat \& J. Tits, {\it Schémas en groupes et immeubles des groupes classiques sur un corps local}, Bull. S.M.F. 112 (1984), 259-301 (= J. Tits, C.P., vol.IV, n°128, 1-43).

   [Ge 81] P. Gérardin, {\it Immeubles des groupes linéaires généraux}, Lect. Notes Math. 880 (1981), 138-178.
   
   [GI 63] O. Goldman \& N. Iwahori, {\it The space of ${\frak p}$-adic norms}, Acta Math. 109 (1963), 137-177.
     
   [La 04] T-Y. Lam, {\it Introduction to Quadratic Forms over Fields}, A.M.S., GSM 67, 2004.
        
   [Sp 55] T.A. Springer, {\it Quadratic forms over fields with a discrete valuation}, Indag. Math. 17 (1955), 352-362 and 18 (1956), 238-246.
   
   [Th 86] J.G. Thompson, {\it Bilinear forms in characteristic $p$ and the Frobenius-Schur indicator}, second part of {\it
   Some finite groups which appear as $\Gal L/K,$ where  $K \subset \Q(\mu$$_n)$}, in {\it Group Theory}, Lect. Notes Math.  1185 (1986), 221-230. 
   
   \bigskip

   Jean-Pierre Serre
   
   Collège de France
   
   3 rue d'Ulm
   
   75005 Paris, France
   
   {\it E-mail}: {\bf jpserre691@gmail.com} 
   
  \end{document}